\documentclass[12pt,a4paper,reqno]{amsart}
\usepackage{lmodern}
\usepackage{amsmath,amsfonts,amssymb,mathrsfs,charter,amsthm}

\newtheorem{theorem}{Theorem}

\newtheorem{corollary}[theorem]{Corollary}

\newcommand{\M}{\mathbb{M}}

\newcommand{\p}{\mathbb{P}}
\newcommand{\tr}{\mathop{{\rm tr}}}
\newcommand{\weaklog}{\mathop{{\rm wlog}}}
\newcommand{\h}{\mathop{{\rm Re}}}

\title{Norm inequalities related to the matrix geometric mean}
\author[Bhatia]{Rajendra Bhatia}
\author[Grover]{Priyanka Grover}
\address{Indian Statistical Institute, Delhi Centre, 7, S.J.S. Sansanwal Marg, New Delhi-110016, India}
\email{rbh@isid.ac.in, pgrover8r@isid.ac.in}

\begin{document}
{\setlength{\baselineskip}%
{1.5\baselineskip}
\begin{abstract}
Inequalities for norms of different versions of the geometric mean of two positive definite matrices are presented.
\end{abstract}

\subjclass[2010]{15A42, 15A18, 47A64, 47A30}

\keywords{Matrix inequalities, Geometric mean, Binomial mean, Log Euclidean mean, Golden-Thompson inequality, Positive definite matrices}
\maketitle

\section{Introduction}
The geometric mean of positive numbers $a$ and $b$ is the number $\sqrt{ab}$, and it satisfies the equations
\begin{equation}
\sqrt{ab}=e^{\frac{1}{2}\left(\log a+\log b\right)}=\lim_{p\rightarrow 0} \left(\frac{a^p+b^p}{2}\right)^{1/p}.
\end{equation}
The quantity 
\begin{equation}
f(p)=\left(\frac{a^p+b^p}{2}\right)^{1/p}, \quad -\infty<p< \infty,
\end{equation}
is called the \emph{binomial mean}, or the \emph{power mean}, and is an increasing function of $p$ on $\left(-\infty,\infty\right)$.

Replacing $a$ and $b$ by positive definite matrices $A$ and $B$, let
\begin{equation}
F(p)=\left(\frac{A^p+B^p}{2}\right)^{1/p}.\label{defnofF}
\end{equation}

In \cite{bhagwat} Bhagwat and Subramanian showed that
\begin{equation}
\lim_{p\rightarrow 0} F(p)=e^{\frac{1}{2}\left(\log A+\log B\right)}.\label{limitofF}
\end{equation}
They also showed that the matrix function $F(p)$ is monotone with respect to $p$, on the intervals $(-\infty,-1]$ and $[1,\infty)$ but not on $\left(-1,1\right)$. (The order $X\leq Y$ on the space $\p$ of $n\times n$ positive definite matrices is defined to mean $Y-X$ is a positive semidefinite matrix.)

The entity in \eqref{limitofF} is called the ``log Euclidean mean'' of $A$ and $B$. However it has some drawbacks, and the accepted definition of the geometric mean of $A$ and $B$ is 
\begin{equation}
A\#_{1/2} B=A^{1/2} \left(A^{-1/2} B A^{-1/2}\right)^{1/2} A^{1/2}.\label{defnofgm}
\end{equation}

It is of interest to have various comparisons between the quantities in \eqref{defnofF}, \eqref{limitofF} and \eqref{defnofgm}, and that is the question discussed in this note.

Generalising \eqref{defnofgm} various authors have considered for $0\leq t\leq 1$
\begin{equation}
A\#_{t}B=A^{1/2}\left(A^{-1/2} B A^{-1/2}\right)^t A^{1/2},\label{geodesic}
\end{equation}
and called it \emph{$t$-geometric mean}, or \emph{$t$-power mean}. In recent years there has been added interest in this object because of its connections with Riemannian geometry \cite{bhatia2}. The space $\p$ has a natural Riemannian metric, with respect to which there is a unique geodesic joining any two points $A,B$ of $\p$. This geodesic can be parametrised as \eqref{geodesic}.

The linear path
\begin{equation}
(1-t)A+tB, \quad 0\leq t\leq 1,\label{line}
\end{equation}
is another path in $\p$ joining $A$ and $B$. It is well known \cite[Exercise 6.5.6]{bhatia2} that
\begin{equation}
A\#_{t}B \leq (1-t)A+tB \text{ for all }0\leq t\leq 1.\label{pathcomparison}
\end{equation}
The special case $t=1/2$ of this is the matrix arithmetic-geometric mean inequality, first proved by Ando \cite{ando1}.

For $0\leq t\leq 1$ let
\begin{equation}
F_t(p)=\left((1-t)A^p +t B^p\right)^{1/p}.\label{defnofFt}
\end{equation}
For $t=1/2$ this is the $F$ defined in \eqref{defnofF}.
It follows from the work in \cite{bhagwat} that 
\begin{equation}
\lim_{p\rightarrow 0} F_t(p)=e^{(1-t) \log A+t \log B},\label{limitofFt}
\end{equation}
and that $F_t(p)$ is monotone with respect to $p$ on  $(-\infty,-1]$ and $[1,\infty)$ but not on $\left(-1,1\right)$. We denote by $\lambda_j(X),\ 1\leq j\leq n$, the decreasingly ordered eigenvalues of a Hermitian matrix $X$, and by $|||\cdot|||$ any unitarily invariant norm on the space $\M$ of $n\times n$ matrices. Our first observation is that while the matrix function $F_t(p)$ defined in \eqref{defnofFt} is not monotone on the whole line $\left(-\infty,\infty\right)$, the real functions $\lambda_j(F_t(p))$ are:

\begin{theorem}\label{1}
Given positive definite matrices $A$ and $B$, let $F_t(p)$ be as defined in \eqref{defnofFt}. Then for $1\leq j\leq n$ the function $\lambda_j(F_t(p))$ is an increasing function of $p$ on $\left(-\infty,\infty\right)$.
\end{theorem}

As a corollary $|||F_t(p)|||$ is an increasing function of $p$ on $\left(-\infty,\infty\right)$. In contrast to this, Hiai and Zhan \cite{zhan} have shown that the function $|||\left(A^p+B^p\right)^{1/p}|||$ is \emph{decreasing} on $(0,1]$ (but not necessarily so on $(1,\infty)$). A several variable version of both our Theorem \ref{1} and this result of Hiai and Zhan can be established (see Remark 1).

Combining Theorem \ref{1} with a result of Ando and Hiai \cite{hiai} we obtain a comparison of norms of the means \eqref{defnofF}, \eqref{limitofF}, \eqref{defnofgm}, and their $t$-generalisations:

\begin{corollary}\label{2}
Let $A$ and $B$ be two positive definite matrices. Then for $p>0$
\begin{equation}
|||A\#_{t}B|||\leq |||e^{(1-t)\log A+t\log B}|||\leq |||\left((1-t)A^{p}+t B^p\right)^{1/p}|||.\label{pinequality}
\end{equation}
\end{corollary}

The first inequality in \eqref{pinequality} is proved in \cite{hiai} as a complement to the famous Golden-Thompson inequality: for Hermitian matrices $H, K$ we have $|||e^{H+K}|||\leq |||e^H e^K|||$. Stronger versions of this inequality due to Araki \cite{araki} and Ando-Hiai \cite{hiai} can be used to obtain a refinement of \eqref{pinequality}. We have for $0\leq t\leq 1$
\begin{eqnarray}
|||A\#_{t}B|||&\leq& |||e^{(1-t)\log{A}+t \log{B}}|||\nonumber\\
&\leq& |||(B^{\frac{tp}{2}} A^{(1-t)p} B^{\frac{tp}{2}})^{1/p}|||\nonumber\\
&\leq& |||((1-t)A^p+tB^p)^{1/p}|||.\label{pinequalities}
\end{eqnarray}

We draw special attention to the case $p=1$ for which further refinements are possible.
\begin{theorem}\label{3}
Let $A$ and $B$ be positive definite matrices. Then
\begin{eqnarray}
|||A\#_{t}B|||&\leq& |||e^{(1-t)\log{A}+t \log{B}}|||\nonumber\\
&\leq& |||B^{\frac{t}{2}} A^{1-t} B^{\frac{t}{2}}|||\nonumber\\
&\leq&  \left|\left|\left|\frac{1}{2}\left(A^{1-t} B^{t}+ B^{t} A^{1-t}\right)\right|\right|\right|\nonumber\\
&\leq& |||A^{1-t} B^t|||\nonumber\\
&\leq& |||(1-t)A+tB|||.\label{inequalities}
\end{eqnarray}
\end{theorem}

For convenience we have stated these results as inequalities for unitarily invariant norms. Many of these inequalities have stronger versions (with log majorisations instead of weak majorisations). This is explained along with the proofs in Section 2. For the special case $t=1/2$ we provide an alternative special proof for a part of Theorem \ref{3}, and supplement it with other inequalities. Section 3 contains remarks and comparisons with known results, some of which are very recent.

\section{Proofs} \emph{Proof of Theorem \ref{1}}\hspace{0.4cm} 
Let $0<p<p'$. Then the map $f(t)= t^{p/p'}$ on $[0,\infty)$ is matrix concave; see \cite[Chapter V]{bhatia1}. Hence 
\begin{equation*}
(1-t)A^p+t B^p\leq \left((1-t)A^{p'}+t B^{p'}\right)^{p/p'}.
\end{equation*}
This implies that
\begin{equation*}
\lambda_j\left((1-t)A^p+t B^p\right)\leq \lambda_j\left((1-t)A^{p'}+t B^{p'}\right)^{p/p'}.
\end{equation*}
Taking $p$th roots of both sides, we obtain
\begin{equation}
\lambda_j\left((1-t)A^p+t B^p\right)^{1/p}\leq \lambda_j\left((1-t)A^{p'}+t B^{p'}\right)^{1/p'}.\label{pthroot}
\end{equation}
Next consider the case $p<p'<0$. Then $0<p'/p<1$.
Arguing as above we obtain
\begin{equation*}
\lambda_j\left((1-t)A^{p'}+t B^{p'}\right)\leq \lambda_j\left((1-t)A^{p}+t B^{p}\right)^{p'/p}.
\end{equation*}
Take $p'$th roots of both sides. Since $p'<0$, the inequality is reversed and we get the inequality \eqref{pthroot} in this case too. Now let $p$ be any positive real number. Using the matrix convexity of the function $f(t)=t^{-1}$ we see that
\begin{equation*}
\left((1-t)A^{-p}+t B^{-p}\right)^{-1}\leq (1-t)A^p+tB^p.
\end{equation*}
From this we get an inequality for the $j$th eigenvalues, and then for their $p$th roots; i.e.,
\begin{equation*}
\lambda_j\left((1-t)A^{-p}+t B^{-p}\right)^{-1/p}\leq \lambda_j\left((1-t)A^p+tB^p\right)^{1/p}.
\end{equation*}
It follows from the above cases that for any $p<0<p'$
\begin{equation}
\lambda_j\left((1-t)A^p+tB^p\right)^{1/p} \leq \lambda_j\left((1-t)A^{p'}+tB^{p'}\right)^{1/p'}.\label{pandp'}
\end{equation}
Taking limit as $p'\rightarrow 0$ and using \eqref{limitofFt} we get
$$\lambda_j\left((1-t)A^p+tB^p\right)^{1/p}\leq \lambda_j\left(e^{(1-t)\log A+t \log B}\right)$$
i.e., for any $p<0$ we have $\lambda_j\left(F(p)\right)\leq \lambda_j\left(F(0)\right).$
For the case $p>0$ a similar argument shows that $\lambda_j\left(F(p)\right)\geq \lambda_j\left(F(0)\right).$ \qed

\bigskip

\emph{Proof of Theorem \ref{3}}\hspace{0.4cm}
The first inequality in \eqref{inequalities} follows from a more general result of Ando and Hiai \cite{hiai}. They showed that for Hermitian matrices $H$ and $K$,
$|||\left(e^{pH}\#_t\ e^{pK}\right)^{1/p}|||$ increases to $|||e^{(1-t)H+tK}|||$ as $p\downarrow 0$. Choosing $H=\log A, \ K=\log B$, and $p=1$, we obtain the first inequality in \eqref{inequalities}. The Golden-Thompson inequality generalised to all unitarily invariant norms (see \cite[p. 261]{bhatia1} says that $|||e^{H+K}|||\leq |||e^{K/2} e^H e^{K/2}|||$. Using this we obtain the second inequality in \eqref{inequalities}. (We remark here that it was shown in \cite{bhatiaemi} that the generalised Golden Thompson inequality follows from a generalised exponential metric increasing property. The latter is related to the metric geometry of the manifold $\p$. So its use in the present context seems natural.) Given a matrix $X$ we denote by $\h X$ the matrix $\frac{1}{2}(X+X^*)$. By Proposition IX.1.2 in \cite{bhatia1} if a product $XY$ is Hermitian, then $|||XY|||\leq |||\h (YX)|||$. Using this we obtain the third inequality in \eqref{inequalities}. The fourth inequality follows from the general fact $|||\h X|||\leq |||X|||$ for all $X$. The last inequality in \eqref{inequalities} is a consequence of the matrix Young inequality proved by T. Ando \cite{ando}.\qed 

For Hermitian matrices $H,K$ let $\lambda_1(H)\geq \cdots\geq \lambda_n(H)$ and $\lambda_1(K)\geq \cdots\geq \lambda_n(K)$ be the eigenvalues of $H$ and $K$ respectively. Then the \emph{weak majorisation} $\lambda(H)\prec_w \lambda(K)$ means that
$$\sum_{i=1}^k \lambda_i(H)\leq \sum_{i=1}^k \lambda_i(K),\quad k=1,2,\ldots,n.$$
If in addition for $k=n$ there is equality here, then we say $\lambda(H)\prec \lambda(K)$.
For $A,B\geq 0$ we write \begin{equation}
\lambda(A)\prec_{\log}\lambda(B)$$ if 
$$\prod_{i=1}^k \lambda_i(A)\leq \prod_{i=1}^k \lambda_i(B), \quad k=1,\ldots,n-1 \label{logmajorisation}
\end{equation}
and $$\prod_{i=1}^n \lambda_i(A)=\prod_{i=1}^n \lambda_i(B), \text{ that is } \det A=\det B.$$
We refer to it as \emph{log majorisation}. We say $A$ is \emph{weakly log majorised} by $B$, in symbols $\lambda(A)\prec_{\weaklog} \lambda(B)$, if \eqref{logmajorisation} is fulfilled. It is known that 
$$\lambda(A)\prec_{\weaklog}\lambda(B) \text{ implies } \lambda(A)\prec_w \lambda(B),$$
so that $|||A|||\leq |||B|||$ for any unitarily invariant norm. (See \cite{bhatia1} for facts on majorisation used here.)

There are stronger versions of some of the inequalities in \eqref{pinequalities}. We have for $p>0$
\begin{align}
\lambda(A\#_t B) &\prec_{\log} \lambda(e^{(1-t)\log A+t\log B})\nonumber\\
&\prec_{\log} \lambda\left(B^{tp/2} A^{(1-t)p} B^{tp/2}\right)^{1/p}\nonumber\\
&= \lambda\left(A^{(1-t)p} B^{tp}\right)^{1/p}\nonumber\\
&\prec_{\weaklog} \lambda\left((1-t) A^p+t B^p\right)^{1/p}.\label{stronginequalities}
\end{align}
The first inequality is a result by Ando and Hiai \cite{hiai}. The second inequality follows from a result by Araki \cite{araki}. The last inequality above follows from the matrix version of Young's inequality by Ando \cite{ando}.

A further strengthening of the first inequality in \eqref{stronginequalities} replacing log majorisation by pointwise domination is not possible. For $t=1/2$ this would have said
$$\lambda_j(A\#_{1/2} B) \leq \lambda_j\left(e^{\frac{\log A+\log B}{2}}\right).$$
This is refuted by the example $A=\left[\begin{array}{cc}
2 & \ 0\\
\ 0 & 1
\end{array}\right]$, $B=\left[\begin{array}{cc}
3 & 3 \\
3 & 9/2
\end{array}\right]$. A calculation shows that $\lambda_2(A\#_{1/2} B)=1$ and $\lambda_2(e^{\frac{\log A+\log B}{2}})\approx 0.9806$.

The case $t=1/2, \ p=1$ is special. Following an idea of Lee \cite{lee} we present a different proof of the majorisation
\begin{equation}
\lambda\left(A\#_{1/2}B\right)\prec_{\log} \lambda\left(B^{1/4}A^{1/2}B^{1/4}\right).\label{majorized}
\end{equation}
The geometric mean $A\#_{1/2} B$ satisfies the equation $A\#_{1/2}B=A^{1/2} U B^{1/2}$ for some unitary U. See \cite[p.109]{bhatia1}. Therefore for the operator norm $\|\cdot\|$ we have
\begin{eqnarray}
\|A\#_{1/2} B\|&=& \|A^{1/2} U B^{1/2}\|\nonumber\\
&=& \|A^{1/4} A^{1/4} U B^{1/4} B^{1/4}\|\nonumber\\
&\leq& \|A^{1/4} U B^{1/4} B^{1/4} A^{1/4}\|\nonumber\\
&\leq& \|A^{1/4} U B^{1/4}\| \|B^{1/4} A^{1/4}\|.\label{inequalities2}
\end{eqnarray}
Here the first inequality is a consequence of the fact that if $XY$ is Hermitian, then $\|XY\|\leq\|YX\|$. Next note that
\begin{eqnarray}
\|A^{1/4} U B^{1/4}\|^2&=& \|A^{1/4} U B^{1/2} U^* A^{1/4}\|\nonumber\\
&\leq& \|A^{1/2} U B^{1/2} U^*\|\nonumber\\
&=& \|A^{1/2} U B^{1/2}\|=\|A\#_{1/2}B\|.\label{inequalities3}
\end{eqnarray}
Again, to derive the inequality above we have used the fact that $\|XY\|\leq \|YX\|$ if $XY$ is Hermitian. From \eqref{inequalities2} and \eqref{inequalities3} we see that 
\begin{equation*}
\|A\#_{1/2}B\|^{1/2} \leq \|B^{1/4} A^{1/4}\|,
\end{equation*}
and hence
\begin{equation}
\|A\#_{1/2} B\|\leq \|B^{1/4} A^{1/4}\|^2=\|B^{1/4} A^{1/2} B^{1/4}\|.\label{strongaujla}
\end{equation}
This is the same as saying that
\begin{equation}
\lambda_1\left(A\#_{1/2} B\right)\leq \lambda_1\left(B^{1/4} A^{1/2} B^{1/4}\right).\label{lambda1}
\end{equation}
If $\wedge^k (X),\ 1\leq k\leq n$, denotes the $k$th antisymmetric tensor power of $X$, then
\begin{equation*}
\wedge^k\left(A\#_{1/2} B\right)=\wedge^k (A) \#_{1/2} \wedge^k(B).
\end{equation*}
So from \eqref{lambda1} we obtain
\begin{equation*}
\lambda_1\left(\wedge^k\left(A\#_{1/2} B\right)\right)\leq \lambda_1\left(\wedge^k(B)^{1/4} \wedge^k(A)^{1/2} \wedge^k(B)^{1/4}\right).
\end{equation*}
This is the same as saying 
\begin{equation}
\prod_{j=1}^k \lambda_j\left(A\#_{1/2} B\right)\leq \prod_{j=1}^k \lambda_j\left(B^{1/4}A^{1/2}B^{1/4}\right),\quad 1\leq k\leq n.\label{productlambda}
\end{equation}
For $k=n$ there is equality here because
\begin{equation*}
\det\left(A\#_{1/2} B\right)=\det\left(A^{1/2} B^{1/2}\right).
\end{equation*}
From \eqref{productlambda} we have the corollary 
\begin{equation}
\lambda\left(A\#_{1/2} B\right)\prec_w \lambda\left(B^{1/4} A^{1/2} B^{1/4}\right).\label{lambda}
\end{equation}
Included in this is the trace inequality
\begin{equation*}
\tr\left(A\#_{1/2} B\right)\leq \tr A^{1/2} B^{1/2}.
\end{equation*}
This has been noted in \cite{lee}. 

\section{Remarks}

\begin{enumerate}

\item[1.] 
Let $A_1,\ldots, A_m$ be positive definite matrices and let $\alpha_1,\ldots,\alpha_m \geq 0$ be such that $\sum \alpha_j=1$. Let
\begin{equation}
F(p)=\left(\alpha_1 A_1^p+\cdots+\alpha_m A_m^p\right)^{1/p}\label{defnofFsevvar}.
\end{equation}
Then by the same argument as in the proof of Theorem \ref{1}, $\lambda_j(F(p))$ is increasing in $p$ on $(-\infty,\infty)$. In particular for $\alpha_1=\cdots=\alpha_m=1/m$ the function $\lambda_j\left(\left(\frac{A_1^p+\cdots+A_m^p}{m}\right)^{1/p}\right)$ is an increasing function of $p$ on $(-\infty,\infty)$. Therefore \\ $\left|\left|\left|\left(\frac{A_1^p+\cdots+A_m^p}{m}\right)^{1/p}\right|\right|\right|$ is an increasing function of $p$ on $(-\infty,\infty)$.
In contrast, it can be shown that $\left|\left|\left|\left(A_1^p+\cdots+A_m^p\right)^{1/p}\right|\right|\right|$ is a decreasing function of $p$ on $(0,1]$. For $m=2$ Hiai and Zhan have shown this using the following result of Ando and Zhan \cite{andozhan}. For positive operators $A,B$ and $r\geq 1$
$$|||(A+B)^r|||\geq |||A^r+B^r||| .$$
A several variable version of this follows from \cite[Theorem 5 (ii)]{kittaneh1} of Bhatia and Kittaneh:
$$|||(A_1+\cdots+A_m)^r|||\geq |||A_1^r+\cdots+A_m^r||| \text{ for } r\geq 1.$$ 
By imitating the argument in \cite{zhan} one can show $\left|\left|\left|\left(A_1^p+\cdots+A_m^p\right)^{1/p}\right|\right|\right|$ is a decreasing function of $p$ on $(0,1]$.

\item[2.] In \cite{bhagwat} Bhagwat and Subramanian showed that for positive definite matrices $A_1,\ldots, A_m$ and $\alpha_1,\ldots,\alpha_m \geq 0$ such that $\sum \alpha_j=1$ 
$$\lim_{p\rightarrow 0} \left(\alpha_1 A_1^p+\cdots+\alpha_m A_m^p\right)^{1/p}=e^{\alpha_1 \log A_1+\cdots+\alpha_m \log A_m}.$$ It follows from Remark 1 that
$$|||e^{\alpha_1 \log A_1+\cdots+\alpha_m \log A_m}|||\leq |||\left(\alpha_1 A_1^p+\cdots+\alpha_m A_m^p\right)^{1/p}||| \text{ for }p>0.$$ 

\item[3.] Recently several versions of geometric mean for more than two positive definite matrices have been considered by various authors. (See \cite{andolimathias}, \cite{bhatiaholbrook} and \cite{bini}.) For positive definite matrices $A_1,\ldots, A_m$ let $G(A_1,\ldots,A_m)$ denote any of these geometric means. Our discussion in Corollary 2 and Remark 2 raises the question whether
$$|||G(A_1,\ldots,A_m)|||\leq |||e^{\frac{\log A_1+\cdots+\log A_m}{m}}|||.$$

\item[4.] By Ando's characterisation of the geometric mean if $X$ is a \emph{Hermitian} matrix and 
$$\left[\begin{array}{ccc} 
A & X\\
X & B
\end{array}\right]\geq 0, \text{ then } X\leq A\# B.$$
Since $$\left[\begin{array}{ccc} 
A & -X\\
-X & B
\end{array}\right]=\left[\begin{array}{ccc} 
I & \ 0\\
\ 0 & -I
\end{array}\right] \left[\begin{array}{ccc} 
A & X\\
X & B
\end{array}\right] \left[\begin{array}{ccc} 
I & \ 0\\
\ 0 & -I
\end{array}\right]$$
 we have
$$\left[\begin{array}{ccc} 
A & -X\\
-X & B
\end{array}\right]\geq 0 \quad \text{ if } \quad \left[\begin{array}{ccc} 
A & X\\
X & B
\end{array}\right]\geq 0.$$
Hence $\pm X \leq A\# B$. Then by \cite[Lemma 2.1]{kittaneh2}, $|||X|||\leq |||A\# B|||$. In contrast to this, we do have that $$\left[\begin{array}{ccc} A & A^{1/2} B^{1/2} \\ B^{1/2} A^{1/2} & B\end{array}\right]=\left[\begin{array}{ccc} A^{1/2} & 0\\ B^{1/2} & 0\end{array}\right] \left[\begin{array}{ccc} A^{1/2} & B^{1/2} \\ 0 & 0\end{array}\right]\geq 0$$ but we have the opposite inequality $|||A\# B|||\leq |||A^{1/2} B^{1/2}|||$. 

\item[5.] Among the several matrix versions of the arithmetic-geometric mean inequality proved by Bhatia and Kittaneh \cite{kittaneh} one says that $4|||AB|||\leq |||(A+B)^2|||$. Using this and Theorem \ref{1} we have
\begin{equation}
|||A^{1/2} B^{1/2}|||\leq \left|\left|\left|\left(\frac{A^{1/2}+B^{1/2}}{2}\right)^2\right|\right|\right|\leq \left|\left|\left|\left(\frac{A^p+B^p}{2}\right)^{1/p}\right|\right|\right| \text{ for } p\geq 1/2.
\end{equation}
For $t=1/2$ this extends the chain of inequalities \eqref{inequalities} in another direction.

\item[6.] In a recent paper \cite{aujla} Matharu and Aujla have shown that
\begin{equation}
\lambda\left(A\#_tB\right)\prec_{\log} \lambda\left(A^{1-t} B^t\right).\label{gmandgeodesic}
\end{equation}
For their proof they use the Furuta inequality. The inequality \eqref{majorized} follows from this. As a corollary these authors observe that
\begin{equation}
|||A\#_{1/2}B|||\leq |||(B^{1/2} A B^{1/2})^{1/2}|||. \label{coraujla}
\end{equation}
In fact, from \eqref{gmandgeodesic} one can deduce the stronger inequality \eqref{lambda}. By IX.2.10 in \cite{bhatia1} we have for $A,B$ positive definite and $0\leq t\leq 1$
\begin{equation*}
|||B^t A^t B^t|||\leq |||(BAB)^t|||.
\end{equation*} 
So, the inequality \eqref{lambda} is stronger than \eqref{coraujla}. In turn, the latter inequality is stronger than one proved by T. Kosem \cite{kosem} who showed 
\begin{equation*}
|||(A\#_{1/2}B)^2|||\leq |||B^{1/2} A B^{1/2}|||. 
\end{equation*}
This follows from \eqref{coraujla} because the majorisation $x\prec_w y$ for positive vectors implies $x^2\prec_w y^2$.

\item[7.] The third inequality in \eqref{inequalities} can be derived from the arithmetic-geometric mean inequality of Bhatia-Davis \cite{davis}
\begin{equation*}
|||A^{1/2}XA^{1/2}|||\leq \left|\left|\left|\frac{1}{2} (AX+XA)\right|\right|\right|
\end{equation*} 
valid for all $X$ and positive definite $A$. There are several refinements of this inequality, some of which involve different means (Heinz means, logarithmic means, etc.) Each such result can be used to further refine \eqref{inequalities}.
\end{enumerate}

\end{document}